\documentclass[letterpaper, 10 pt, conference]{ieeeconf}  % Comment this line out
\IEEEoverridecommandlockouts                              % This command is only
\usepackage{graphics} % for pdf, bitmapped graphics files
\usepackage{epsfig} % for postscript graphics files
\usepackage{times} % assumes new font selection scheme installed
\usepackage{amsmath} % assumes amsmath package installed
\usepackage{amssymb}  % assumes amsmath package installed
\usepackage{bbm}
\usepackage{color}

\newtheorem{lem}{Lemma}
\newtheorem{cor}{Corollary}

\newtheorem{rem}{Remark}

\title{\LARGE \bf
Maximum Throughput Problem in Dissipative Flow Networks with Application to Natural Gas Systems
}

\author{Sidhant Misra, Marc Vuffray, Michael Chertkov% <-this % stops a space
\thanks{S. Misra, M. Vuffray and M. Chertkov. are with the Center for Nonlinear Studies and Theoretical Division T-4 of Los Alamos National Laboratory,
        Los Alamos, NM 87544. M. Chertkov is also affiliated with the New Mexico Consortium, Los Alamos, NM 87544.
        {\tt\small \{sidhant|vuffray|chertkov\}@lanl.gov}}
}

\begin{document}

\maketitle
\thispagestyle{empty}
\pagestyle{empty}

%%%%%%%%%%%%%%%%%%%%%%%%%%%%%%%%%%%%%%%%%%%%%%%%%%%%%%%%%%%%%%%%%%%%%%%%%%%%%%%%
\begin{abstract}
We consider a dissipative flow network that obeys the standard linear nodal flow conservation, and where flows on edges are driven by potential difference between adjacent nodes. We show that in the case when the flow is a monotonically increasing function of the potential difference, solution of the network flow equations is unique and can be equivalently recast as the solution of a strictly convex optimization problem. We also analyze the maximum throughput problem on such networks seeking  to maximize the amount of flow that can be delivered to the loads while satisfying bounds on the node potentials. When the dissipation function is differentiable we develop a representation of the maximum throughput problem in the form of a twice differentiable biconvex optimization problem exploiting the variational representation of the network flow equations.  In the process we prove a special case of a certain monotonicity property of dissipative flow networks. When the dissipation function follows a power law with exponent greater than one, we suggest a mixed integer convex relaxation of the maximum throughput problem.  Finally, we illustrate application of these general results to balanced, i.e. steady, natural gas networks also validating the theory results through simulations on a test case.
\end{abstract}
%%%%%%%%%%%%%%%%%%%%%%%%%%%%%%%%%%%%%%%%%%%%%%%%%%%%%%%%%%%%%%%%%%%%%%%%%%%%%%%%
\begin{section}{Introduction}
We consider a network on  an undirected graph $\mathcal{G} = (\mathcal{V}, \mathcal{E})$, where $\mathcal{V}$  is the set of nodes with $|\mathcal{V}| = m$ and $\mathcal{E}$ is the set of oriented edges with $| \mathcal{E}|  = n$. For a dissipative flow network, we have two sets of variables describing the network. Each $(i,j) \in \mathcal{E}$ is assigned a nominal direction $i \rightarrow j$ and is characterized by a flow $\phi_{ij} \in \mathbbm{R}$ from node 
$i$ to node $j$, and each node $i \in \mathcal{V}$ is characterized by a node potential $\pi_i \in \mathbbm{R}$. Such networks which have a description via potentials are determined by a set of two equations:``flow conservation" and ``potential drop relation." Perhaps the most famous example of the dissipative flow networks is a linear DC resistive electric circuit, where the flows $\phi_{ij}$ are currents and node potentials $\pi_i$ are voltage potentials. In this case, the flow conservation equations correspond to the so-called Kirchoff's first law, and the potential drop relation is the Kirchoff's second law, or Ohm's law. Resistive circuits are linear dissipative networks, where
the potential drop relation given by the Ohm's law is linear. In this paper we consider more general dissipative networks where the potential drop relations is a general nonlinear relation described by a monotonically increasing dissipation function.

We formulate a variational representation of the network flow equations in the form of a  strictly convex optimization problem, thus proving the uniqueness of the solution. Additionally, this solution can be obtained efficiently by using convex optimization algorithms such as the interior point method. The energy function method for gas flow networks was studied in \cite{77Nau,12BNV}. The maximum throughput problem is an optimization problem that aims at maximizing the amount of flow delivered from sources to loads, while satisfying bounds on the potentials and possibly other control variables. For general dissipative flow networks, this problem is non-convex. In the case when the dissipation function is differentiable, we develop a twice differentiable optimization formulation of the maximum throughput problem that is biconvex in the optimization variables in a certain regime. We obtain several properties of the primal and dual form
of the energy function that are useful to provide gradient and hessian information to optimization software during numerical implementation. In the process we prove a special case of the so-called ``monotonicity property" of dissipative flow networks, also discussed in a companion paper \cite{Vuffray15}.

When the dissipation function takes the form of a power law with exponent greater than one, we provide a mixed integer convex relaxation of the maximum throughput problem, where the degree of non-linearity is governed by the exponent. Together, we have a framework to obtain a feasible solution using the energy function formulation, and then check the quality of the solution (degree of sub-optimality) using the mixed integer convex formulation. As an example, we demonstrate how the setting applies to the natural gas networks and also illustrate its performance on a numerical example.

The rest of the manuscript is organized as follows. In Section \ref{sec:dissipative} we define dissipative flow networks and equations that govern them. Section \ref{max_throughput} states the maximum throughput problem for a dissipative flow network. Section \ref{energy_function_optimization} introduces variational formulations of the network flow equations and their properties and provides a reformulation of the maximum throughput problem using the energy functions. In Section \ref{MICP} we introduce a mixed integer convex relaxation for the maximum throughput problem where the dissipation function follows a power law with exponent greater than one. In Section
\ref{ng_networks} we briefly describe natural gas networks and adopt the techniques of the previous Sections to this special enabling case. Finally, in Section \ref{simulations} we provide computational results on natural gas network models and then summarise and present path forward in Section \ref{sec:conclusions}.

\end{section}
%%%%%%%%%%%%%%%%%%%%%%%%%%%%%%%%%%%%%%%%%%%%%%%%%%%%%%%%%%%%%%%%%%%%%%%%%%%%%%%%
\begin{section}{Dissipative Flow Networks} \label{sec:dissipative}
A dissipative flow network, where flows on edges are driven by potential difference between adjacent nodes, is described by the following two sets of the Network Flow (NF) equations:
\begin{itemize}
	\item[(a)] Flow conservation:
		\begin{align}
			\sum_{j \in \partial i} \phi_{ij} = q_i, \ \forall i \in \mathcal{V} \label{flow_conservation}
		\end{align}
		where $q_i \in \mathbbm{R}$ is the injection at node $i$.
	\item[(b)] Potential drop relation:
		\begin{align}
			\pi_i - \pi_j = f_{ij}(\phi_{ij}) - b_{ij}, \  \forall (i,j) \in \mathcal{E}. \label{potential_flow_relation}
		\end{align}
		The function $f_{ij}(.)$ is a differentiable, strictly monotonically increasing odd function which is called the \textit{dissipation function}.
\end{itemize}
\vspace{0.1in}
Several networks can be described within the framework. For example, DC resistive electric circuits \cite{98Bol}, AC power flows in the lossless approximation \cite{58Ayl,81ASV,81BH,15DLS}, balanced natural gas networks \cite{Osiadacz_1987,Wu2000197,12BNV,Misra2015}, traffic flows \cite{10CSADF,13Var_a,13Var_b} etc. The term $b_{ij}$ allows us to model additional network elements in the actual networks such as compressors in gas networks, voltage transformers in DC resistive circuits  and phase shifters in the approximate AC power flows. %\textbf{unsure}

The system of potential drop relations (\ref{potential_flow_relation}) is translation invariant, i.e., if $(\pi_{i}, \ i \in \mathcal{V})$ is a solution of the NF Eqs.~(\ref{flow_conservation},\ref{potential_flow_relation}), then a translation of the potentials by a constant, $\tilde{\pi}_i = \pi_i + c$, is also a solution. To remove this degree of freedom, the potential at a particular node (node $1$) is fixed at a certain value $\pi_1 = c$. In the case of the resistive circuits, this node is called ``ground", in AC power system respective node is called ``slack bus", and this term is also adopted to the case of the gas networks,  where the slack bus (where the pressure is fixed) is usually assigned to the largest source of gas in the system. Assuming that parameters $b_{ij}$ are known, we are now left with the number of variables, ($|\mathcal{V}| + |\mathcal{E}| - 1$), equal to the number of equations.  However, this still does not guarantee existence and uniqueness of the NF Eqs.~(\ref{flow_conservation},\ref{potential_flow_relation}).
But as we will see shortly, in Section \ref{sec:primal_energy_function}, this is true for dissipative networks with a strictly monotonically increasing dissipation function. Further, observe that the set of flow conservation equations (\ref{flow_conservation}) automatically implies total flow balance
\begin{align}
	\sum_{i \in \mathcal{V}} q_i = 0. \label{total_flow_balance}
\end{align}

\subsection{Primal Energy Function and Uniqueness of the Network Flow Solution} \label{sec:primal_energy_function}

Let $F_{ij}(x)$ be the anti-derivative of $f_{ij}(x)$, i.e, $\frac{d}{dx} F_{ij}(x) = f_{ij}(x)$.
The following useful remark follows immediately from the fact that $f_{ij}(x)$ is monotonically increasing.
\vspace{0.1in}
\begin{rem} \label{convexity}
	The function $F_{ij}(x)$ is convex.
\end{rem}
\vspace{0.1in}

The set of NF equations (\ref{flow_conservation}) and (\ref{potential_flow_relation}) can be rewritten in the variational form as the optimal solution of the following convex optimization problem \cite{77Nau,12BNV}:
\begin{align}
	\min \quad & \sum_{(i,j) \in \mathcal{E}} F_{ij}(\phi_{ij}) - b_{ij} \phi_{ij}, \label{primal_objective} \\
		\mbox{s.t.} \quad & \sum_{j \in \partial i} \phi_{ij} = q_i. \label{primal_constraint}
\end{align}
The optimization problem in (\ref{primal_objective})-(\ref{primal_constraint}) is clearly convex, since  by Remark \ref{convexity} the cost function is convex, and the constraints are linear.
To see why the solution to (\ref{primal_objective})-(\ref{primal_constraint}) yields the NF equations, we can first form the Lagrangian as
\begin{align}
	L(\phi, \pi) \!=\! \sum_{(i,j) \in \mathcal{E}} F_{ij}(\phi_{ij}) \!-\! b_{ij} \phi_{ij} \!-\! \sum_{i \in \mathcal{V}}  \pi_i\left( \sum_{j \in \partial i} \phi_{ij} - q_i  \right),
\end{align}
where $\pi_i$ are the dual variables corresponding to the constraints (\ref{primal_constraint}). To minimize the Lagrangian, we take its derivative w.r.t. $\phi$ to zero to get
\begin{align}
	f_{ij}(\phi_{ij}) = \pi_i - \pi_j + b_{ij}, \label{lagrangian_optimum}
\end{align}
which means that the dual variables serve the role of the node potentials. Since the optimization problem in (\ref{primal_objective})-(\ref{primal_constraint}) is strictly convex, we also conclude that the NF Eqs.~(\ref{flow_conservation}) and (\ref{potential_flow_relation}) have a unique solution.

\subsection{Dual Energy Function Representation}

We now present another variational formulation of the NF equations via the dual of the optimization problem in (\ref{primal_objective})-(\ref{primal_constraint}). Since $f_{ij}(x)$ is monotonically increasing, it is invertible. Let $g_{ij}(y) \triangleq f_{ij}^{-1}(y)$ denote the inverse of $f_{ij}(x)$.
Let $G_{ij}(y) = F(f_{ij}^{-1}(y)) - y f_{ij}^{-1}(y)$ be the anti-derivative of $g_{ij}(y)$.
The dual of (\ref{primal_objective})-(\ref{primal_constraint}) can be expressed in terms of the node potential variables $\pi_i$ as:
\begin{align}
	\max \sum_{i \in \mathcal{V}} \pi_i q_i - \sum_{(i,j) \in \mathcal{E}}   G_{ij}(\pi_i - \pi_j + b_{ij}). \label{dual_objective}
\end{align}
The above problem, when stated as a minimization, is a convex optimization problem.
Since the Slater's condition \cite{BoydCO} always holds for an underdetermined system of linear equality constraints, we have strong duality and the optimal value of (\ref{primal_objective})-(\ref{primal_constraint}) and (\ref{dual_objective})
are the same.
\end{section}

%%%%%%%%%%%%%%%%%%%%%%%%%%%%%%%%%%%%%%%%%%%%%%%%%%%%%%%%%%%%%%%%%%%%%%%%%%%%%%%%
\begin{section}{The Maximum Throughput Problem} \label{max_throughput}
In a dissipative network the nodes with $q_i > 0$ are called sources and those with $q_i \leq 0$ are called sinks.
Assume that there are constraints on the minimum and maximum value of the allowed node potentials, $\pi_i$, in the form
\begin{align}
	\underline{\pi}_i \leq \pi_i \leq \overline{\pi}_i. \ \forall i \in \mathcal{V}. \label{box_constraints}
\end{align}
The maximum throughput problem seeks to maximize the amount of flow that is delivered to the loads while still satisfying the box constraints on node potentials in (\ref{box_constraints}). Formally:

 \vspace{0.1in} \emph{Maximum Throughput in Dissipative Flow Networks:}
 \begin{align}
 	\min \quad &\sum_{i \in \mathcal{V}} c_i x_i  \label{max_throughput_objective} \\
	\mbox{s.t.} \quad & \sum_{j \in \partial i} \phi_{ij} = q_i, \ \forall i \in \mathcal{V}, \\
		& \pi_i - \pi_j = f_{ij}(\phi_{ij}) - b_{ij}, \  \forall (i,j) \in \mathcal{E},  \label{non_linear} \\
		& \underline{\pi}_i \leq \pi_i \leq \overline{\pi}_i. \ \forall i \in \mathcal{V}. \label{max_throughput_box}
 \end{align}
 The potential drop relations (\ref{non_linear}) are in general non-convex and make the maximum throughput problem a non-convex problem.
\end{section}
%%%%%%%%%%%%%%%%%%%%%%%%%%%%%%%%%%%%%%%%%%%%%%%%%%%%%%%%%%%%%%%%%%%%%%%%%%%%%%%%

\begin{section}{Representing the maximum throughput problem using energy functions } \label{energy_function_optimization}

In this Section, we provide another alternative (but, obviously, equivalent) formulation for the maximum throughput problem using the variational representation of the NF equations in terms of the energy function. To achieve this goal we, first, discuss optimal solution and value at the optima of the energy function.
	
\subsection{Properties of the energy functions}

Define the following convex function given by
	\begin{align}
		E(\pi, b) \triangleq \sum_{(i,j) \in \mathcal{E}} G_{ij}(\pi_i - \pi_j + b_{ij}). \label{E}
	\end{align}
	and let
	\begin{align}
		E^{*}(x, b) \triangleq \sup_{\pi} \left[ \pi^T x - E(\pi, b) \right]. \label{Estar}
	\end{align}	
	The function $E^*(x, b)$ is the convex conjugate of $E(\pi,b)$ with respect to $\pi$. Analogously, define
	\begin{align}
		F(\phi) \triangleq \sum_{(i,j) \in \mathcal{E}} F_{ij}(\phi_{ij}). \label{F}
	\end{align}
	and
	\begin{align}
		F^*(b, q) \triangleq \sup_{\phi}  \quad &\sum_{(i,j) \in \mathcal{E}} b_{ij} \phi_{ij} - F(\phi),  \label{Fstar} \\
				\mbox{s.t.} \quad & \sum_{j \in \partial i} \phi_{ij} = q_i.
	\end{align}
	
The following remark follows directly from the definitions above.
	\begin{rem}
		The function $E^{*}(q, b)$  is equal to the optimal value of the primal-dual pair in  (\ref{primal_objective})-(\ref{primal_constraint}) and (\ref{dual_objective}). Additionally $F^*(b, q) = - E^{*}(q, b)$.
	\end{rem}
	\vspace{0.1in}
	\begin{lem} \label{duality_results}
		The following statements follow from standard results of the convex conjugate and the duality theory \cite{Rockafellar70}.
		\begin{itemize}

			\item[(a)] Let $\phi^*$ be the unique optimum of Eqs.~(\ref{primal_objective})-(\ref{primal_constraint}) and let $\pi^{*}$ be the unique optimum of the Eq.~(\ref{dual_objective}). Then,
					\begin{align}
						\nabla_{x} E^{*}(x, b) &= \pi^{*}, \label{pistar} \\
						\nabla_{b} E^{*}(x,b) &= - \phi^{*}. \label{phistar}
					\end{align}
					
			\item[(b)] The Hessian of $E^*(x, b)$ can be computed with respect to $x$ as follows
				\begin{align}
					\nabla_{xx}^2 E^{*}(x, b)  = \left( \nabla_{\pi \pi}^2E(\pi, b)\right)^{-1},
				\end{align}		
				where
				\begin{align}
					\left[\nabla_{\pi \pi}^2 E(\pi, b)\right]_{ij} &=  - g'(\pi_i - \pi_j + b_{ij}),  \\
					&\quad 1 \leq i \neq j \leq n-1, \ (i,j) \in \mathcal{E}, \label{hessian_off_diag} \\
					\left[\nabla_{\pi \pi}^2 E(\pi, b)\right]_{ii} &=   \sum_{j \in \partial i} g'(\pi_i - \pi_j + b_{ij}), \\
					& \quad \quad  1 \leq i \leq n-1. \label{hessian_diag}
				\end{align}
			
			\item[(c)] The inequality
				\begin{align}
					E(\pi, b) + E^{*}(x, b)  - \pi^T x\leq 0 \label{efunc_inequality}
				\end{align}
				has exactly one solution  given by $\pi^* = \nabla_{x} E^{*}(x, b) $.
		\end{itemize}
	\end{lem}
	\vspace{0.1in}
	
	 \begin{proof}
		\begin{itemize}
			\item[(a)] The statement in Eq.~(\ref{pistar}) follows directly from the theory of convex conjugates. The proof of (\ref{phistar}) is similar, but we include it here for completeness.
					Recall that $E^*(x, b) = -F^*(b,x)$. So, it is enough to prove that $\nabla_{b} F^{*}(b,x) = \phi^*$. By definition of $\phi^*$, we have
					\begin{align*}
						F^*(x, b) = \sum_{(i,j) \in \mathcal{E}} b_{ij} \phi_{ij}^* - F(\phi^*).
					\end{align*}
					So,
					\begin{align*}
						\frac{\partial F^*(x,b)}{\partial b_{kl}} =  \phi_{kl}^* +  \sum_{(i,j) \in \mathcal{E}} b_{ij} \frac{\partial \phi_{ij}^*}{b_{kl}}    - f_{ij}(\phi_{ij}^*) \frac{\partial \phi_{ij}^*}{b_{kl}}
					\end{align*}
					Using Eq.~(\ref{lagrangian_optimum}), one derives
					\begin{align}
						\frac{\partial F^*(x,b)}{\partial b_{kl}}  &=  \phi_{kl}^* +  \sum_{(i,j) \in \mathcal{E}}\frac{\partial \phi_{ij}^*}{b_{kl}} \left( \pi_j  - \pi_i  \right), \\
						&=  \phi_{kl}^* + \sum_{i \in\mathcal{V}}  \pi_i \sum_{j \in \partial i} \frac{\partial \phi_{ij}^*}{b_{kl}}. \label{proof_final_line}
					\end{align}
					Since $\sum_{j \in \partial i} \phi_{ij}^* = x_i$, we have that $\sum_{j \in \partial i} \frac{\partial \phi_{ij}^*}{b_{kl}} = 0$. The proof is complete by substituting this result in Eq.~(\ref{proof_final_line}).
					\vspace{0.05in}
			\item[(b)] Follows from standard properties of convex conjugates and direct computations.
			\item[(c)] By definition $-E^*(x, b)$ is the minimum of $E(\pi, b)$ w.r.t. $\pi$. Since $E(\pi, b)$ is strictly convex, only its unique minima $\pi^*$ can satisfy (\ref{efunc_inequality}).
		\end{itemize}
	\end{proof}
	
	\noindent For the statement in Lemma \ref{duality_results} (b) to be well-defined, we need that the matrix $\nabla^2_{\pi\pi} E(\pi,b)$ is invertible. Note that since the pressure at node $0$ is fixed at $\pi_0 = c$,
	indices in Eq.~(\ref{hessian_off_diag}) and Eq.~(\ref{hessian_diag}) run from $1$ to $n-1$. From the expressions for the elements of the hessian, we derive that $\nabla^2_{\pi \pi}E(\pi, b)$ is a $(n-1)\times(n-1)$
	sub-matrix of a weighted graph {Laplacian} matrix with edge weights given by $w_{ij} = g'(\pi_i - \pi_j + b_{ij})$. Since the graph $\mathcal{G}$ is assumed to be connected, it follows that $\nabla_{\pi \pi}^2 E(\pi, b)$
	has full rank and hence is invertible \cite{Merris94}, \cite{Chaiken78}.
	
Let us now take a slight detour from the main goal (reformulation of the maximum throughput problem) and state a certain ``monotonicity property" that holds for dissipative FN. For a more general statement and proof of the monotonicity property without assuming the differentiability of the
	dissipation function $f_{ij}(\phi_{ij})$ as well as its implication on robust optimization on dissipative flow networks, we refer to \cite{Vuffray15}.
	
	\vspace{0.1in}
	\begin{cor} (Monotonicity Property) \label{monotonicity_property}
		Consider a dissipative flow network where the parameters $b_{ij}$ and the potential $\pi_0$ are fixed.
		 Let $q^{(1)} = (q_1^{(1)}, \ldots, q_{m-1}^{(1)})$ and $q^{(2)} = (q_1^{(2)}, \ldots, q_{n-1}^{(2)})$ be two sets of injections. The injection at node $0$ is automatically determined
		by the total flow balance Eq.~(\ref{total_flow_balance}). Assume that for all $1 \leq i \leq n-1$ one has $q_i \leq r_i$. Let $\pi^{(1)} = (\pi_1^{(1)}, \ldots \pi_{n-1}^{(1)})$ and
		$\pi^{(2)} = (\pi_1^{(2)}, \ldots \pi_{n-1}^{(2)})$ be the corresponding potentials resulting from the two injection profiles. Then one arrives at $\pi_i^{(1)} \leq \pi_i^{(2)}$ for all $1 \leq i \leq n-1$.
	\end{cor}
	\vspace{0.1in}
	
	\subsection{Equivalent formulation of the maximum throughput problem}
	
	\vspace{0.1in}
	Equipped with Lemma \ref{duality_results}, one restates the maximum throughput problem (\ref{max_throughput_objective})-(\ref{max_throughput_box}) as
	\begin{align}
		\min \quad & c^T x,  \label{reform_objective} \\
		\mbox{s.t.} \quad & E(\pi, b) + E^*(x, b) \leq 0, \label{efunc_constraint} \\
			& \underline{\pi}_i \leq \pi_i \leq \overline{\pi}_i. \ \forall i \in \mathcal{V}. \label{reform_box}
	\end{align}
		This formulation is exactly equivalent to the original maximum throughput formulation. Let us present two other alternative formulations,  which (as we see below) may have some practical advantages. \\

	\noindent	\textit{Energy function formulation $1$:}
	\begin{align}
		\min \quad & c^T x,  \label{reform_objective_practical1} \\
		\mbox{s.t.} \quad & E(\pi, b) + E^*(x,b) \leq \epsilon, \nonumber \\
			& \underline{\pi}_i \leq \pi_i \leq \overline{\pi}_i. \ \forall i \in \mathcal{V}. \nonumber
	\end{align}
	\textit{Energy function formulation $2$:}
	\begin{align}
		\min \quad & c^T x + M(E(\pi,b) + E^*(x,b) - \pi^Tx) \label{reform_objective_practical2}   \\
		\mbox{s.t.} \quad& \underline{\pi}_i \leq \pi_i \leq \overline{\pi}_i. \ \forall i \in \mathcal{V}. \nonumber
	\end{align}
	In (\ref{reform_objective_practical1}) $\epsilon > 0$ is a small constant.
	In (\ref{reform_objective_practical2}) $M$ is a large positive constant that helps to enforce the constraint (\ref{efunc_constraint}).
	
Eq.~(\ref{reform_objective_practical1}) and Eq.~(\ref{reform_objective_practical2}) are useful because,
	 first, { the two expressions} are twice differentiable and hence allow application of second order gradient descent methods for local optimization. Second, when  $b$ is not an optimization variable, {the expressions are} biconvex in $\pi$ and $x$ which means that we can use here standard heuristics, such as alternating minimization over $\pi$ and $x$. Moreover, solvers using second order methods can be provided with first and second derivative information by using Eqs.~(\ref{pistar})-(\ref{phistar}) and Eqs.~(\ref{hessian_off_diag})-(\ref{hessian_diag}). The biconvexity is lost when $b$ is also subject to optimization, but the twice differentiability remains.

	\vspace{0.1in}
\end{section}
%%%%%%%%%%%%%%%%%%%%%%%%%%%%%%%%%%%%%%%%%%%%%%%%%%%%%%%%%%%%%%%%%%%%%%%%%%%%%%%%
\begin{section}{Mixed Integer Convex Relaxation}	\label{MICP}

In this Section we formulate a mixed integer convex relaxation of the maximum throughput problem  (\ref{max_throughput_objective})-(\ref{max_throughput_box}) in the special case when the
	dissipation function $f(x)$ is a power law, i.e.,  $f(x) = sgn(x) |x|^{\alpha}$ with $\alpha \geq 1$. Let us start introducing binary variables $s_{ij} \in \{ -1, 1 \}$ associated with each edge $(i,j) \in \mathcal{E}$ that denote the flow direction,
	i.e., $s_{ij} \phi_{ij} \geq 0$. We can then restate the non-convex constraint (\ref{non_linear}) as
	\begin{align}
		&\pi_i - \pi_j  + b_{ij} = s_{ij} |\phi_{ij}|^{\alpha}, \label{binary_flow_variables} \\
		\mbox{or equivaliently } &s_{ij}(\pi_i - \pi_j+ b_{ij}) = |\phi_{ij}|^{\alpha}.
	\end{align}
	Next, we apply the McCormick relaxation \cite{McCormick76} to the bilinear term on the lhs and relax the equality to an inequality to get the following pair of inequality relaxations of (\ref{non_linear}):
	\begin{align}
		|\phi_{ij}|^{\alpha} &\leq s_{ij} (\overline{\pi}_i - \underline{\pi}_j)  - (\pi_{i} - \pi_j) +  (\overline{\pi}_i - \underline{\pi}_j),  \label{MICC1} \\
		|\phi_{ij}|^{\alpha} &\leq s_{ij} (\underline{\pi}_i - \overline{\pi}_j ) + (\pi_i - \pi_j) - (\underline{\pi}_i - \overline{\pi}_j). \label{MICC2}
	\end{align}
	Since $\alpha \geq 1$ the constraints (\ref{MICC1}) and (\ref{MICC2}) are mixed integer convex constraints. The McCormick relaxation is tight when at least one of the variables is binary, the only relaxation in this  set of constraints consists in replacing the equality in (\ref{binary_flow_variables}) by an inequality.
	
The mixed integer convex relaxation along with the formulation in Section \ref{energy_function_optimization} provide the framework for finding a feasible solution of the maximum throughput problem as well as for testing the quality of the solution through comparison with the lower bound.
	
\end{section}
%%%%%%%%%%%%%%%%%%%%%%%%%%%%%%%%%%%%%%%%%%%%%%%%%%%%%%%%%%%%%%%%%%%%%%%%%%%%%%%%

%%%%%%%%%%%%%%%%%%%%%%%%%%%%%%%%%%%%%%%%%%%%%%%%%%%%%%%%%%%%%%%%%%%%%%%%%%%%%%%%
\begin{section}{Natural Gas Networks}  \label{ng_networks}

In this Section, we introduce natural gas networks and their properties and apply the techniques of Section \ref{energy_function_optimization} and \ref{MICP} to the maximum throughput problem over these networks.
	Natural gas has emerged as an alternative to other fuel sources such as coal and fuel oil due to its lower carbon footprint.
	 Improved technology for tapping new natural gas reserves has increased the availability and reduced its cost
	leading to  increasing utilization of natural gas in power generation (gas turbines) and household heating.
	Natural gas networks use a system of high pressure transmission pipelines to transport gas from sources to loads. The flow of gas in pipelines is driven by gradient in pressure along the length of the pipe. Pressure falls rapidly with distance along the pipelines, so compressor stations are installed to locally boost the pressure and maintain the desired pressure and throughput at the load end.
	Typically, compressor stations are placed every 50-100 miles. Many of these compressors are powered by a fraction of the gas in the transmission line itself. The fraction of gas consumed can vary from $2 - 5 \%$.
	Due to the increasing demand for natural gas, there is an increasing need for optimization of pipeline operational planning. In what follows, we describe two such optimization problems of interest.
	
	The first is the problem of minimizing the cost associated with compression. In the literature, study of this problem can be traced back to \cite{68WL} (see also \cite{10Bor} for a review of this problem) where a Dynamic Programming approach was used to solve it.
	More recently in \cite{Misra13}, a convex Geometric Programming (GP) based approach was developed for compressor cost optimization in tree networks. This approach offers several desirable properties such as fast convergence, and scope for generalization to a heuristic for loopy networks. Since the compressors consume a part of the gas itself, the compressor cost optimization problem is also related to the next problem below.
	
	The second  is the maximum throughput problem introduced in Section \ref{max_throughput}. This is particularly relevant in, for example, Norway, which produces more gas than can be domestically consumed. This situation creates strong economic incentives to sell the excess gas to the rest of Europe. In this case the primary objective of the gas pipeline operator becomes to maximize the throughput of the network within the engineering and contractual constraints that must be satisfied.

In the rest of this Section, we describe the equations governing the flow of natural gas over natural gas networks in the balanced (steady state) regime and apply the machinery of Section \ref{energy_function_optimization} and \ref{MICP} to this special case.
	
	Natural gas networks in the steady state are governed by the following set of NF equations:
	\begin{align}
		\sum_{j \in \partial i} \phi_{ij} &= q_i,  \label{gas_flow_conservation} \\
		p_i^2 - p_j^2 + b_{ij} &= \delta_{ij} \phi_{ij} | \phi_{ij}|, \label{pressure_flow_relation}
	\end{align}
	where $\phi_{ij}$ is the flux of gas (flow per unit length) from node $i$ to node $j$ in the steady-state,
	and $p_i$ is the pressure at node $i$. The quantity $b_{ij}$ models the pressure boost provided by the compressor on the edge $(i,j) \in \mathcal{E}$.
	For a more detailed derivation of the steady-state gas flow equations from dynamic fluid mechanic equations, we refer to \cite{Misra13}.

From Eqs.~(\ref{gas_flow_conservation},\ref{pressure_flow_relation}) it is apparent that natural gas networks in steady state can be represented within the framework of dissipative flow networks introduced
	 in Section \ref{sec:dissipative} by setting
	\begin{align}
		\pi_i  &= p_i^2, \\
		\mbox{and }f_{ij}(\phi_{ij}) &= \delta_{ij} \phi_{ij} |\phi_{ij}|.
	\end{align}

	\subsection{Energy function representation for maximum throughput in natural gas networks}	

	In this Section, we adopt the reformulation of the maximum throughput problem from Section \label{energy_function_optimization} to the case of the natural gas networks.
	The expressions for functions $E$ and $F$ in Eqs.~(\ref{E},\ref{F}) for gas networks are given by
	\begin{align}
		E(\pi, b) &= \frac{2}{3} \sum_{(i,j) \in \mathcal{E}} \frac{1}{\sqrt{\delta_{ij}}} | \pi_i - \pi_j + b_{ij}|^{\frac{3}{2}}, \\
		F(\phi) &= \frac{1}{3} \sum_{(i,j) \in \mathcal{E}} \delta_{ij}  | \phi_{ij} |^3.
	\end{align}
	The inverse $g_{ij}(y)  = f_{ij}^{-1}(y)$ is given by
	\begin{align}
		g_{ij}(y) = sgn(y) \sqrt{|y |}
	\end{align}
	and the hessian matrix $\nabla_{\pi \pi}^2 E(\pi, b)$ is given by
	\begin{align}
		\left(\nabla_{\pi \pi}^2 E(\pi, b) \right)_{ij} &= -\frac{1}{2} |\pi_i - \pi_j + b_{ij} |^{-1/2}, \\
		\left(\nabla_{\pi \pi}^2 E(\pi, b) \right)_{ii} &= \sum_{j \in \partial i} \frac{1}{2} |\pi_i - \pi_j + b_{ij} |^{-1/2}.
	\end{align}
	The energy function formulation of the maximum throughput problem for the gas networks can be obtained substituting the expressions above to Eqs.~(\ref{reform_objective}-\ref{reform_box}).
	
	\subsection{Mixed Integer Quadratic Programming (MIQP) relaxation of the maximum throughput in the natural gas networks}

The gas network dissipation function (\ref{pressure_flow_relation}) follows a power law with $\alpha = 2$. So the Mixed Integer Convex relaxation of Section \ref{MICP} is now a Mixed Integer Quadratic Program (MIQP). We write the full program below for completeness.
	
 \vspace{0.2in}
	\textit{MIQP for Maximum Throughput in natural gas networks}
	\begin{align}
		\min \quad & \sum_{i \in \mathcal{V}} c_i x_i  \label{MIQC_cost} \\
		\mbox{s.t.} \quad & \sum_{j \in \partial i} \phi_{ij} = x_i, \\
		&\delta_{ij}  \phi_{ij}^{2} \leq s_{ij} (\overline{\pi}_i - \underline{\pi}_j)  - (\pi_{i} - \pi_j) +  (\overline{\pi}_i - \underline{\pi}_j),  \label{MIQC1} \\
		& \delta_{ij} \phi_{ij}^{2} \leq s_{ij} (\underline{\pi}_i - \overline{\pi}_j ) + (\pi_i - \pi_j) - (\underline{\pi}_i - \overline{\pi}_j), \label{MIQC2} \\
		&\underline{\pi} \leq \pi \leq \overline{\pi}, \\
		& \underline{x} \leq x \leq \overline{x}, \\
		& \underline{b} \leq b \leq \overline{b}. \label{MIQC_end}
	\end{align}
	
As remarked before, the only relaxation in the above program is the conversion of equality into inequality in Eqs.~(\ref{MIQC1},\ref{MIQC2}). This corresponds to allowing ``decompression", i.e., the pressure drop is now
allowed to be larger than that implied by the flow. A similar relaxation was used in the Geometric Programming approach in \cite{Misra13} to preserve convexity, and was observed to perform well in numerical experiments.
	\end{section}

%%%%%%%%%%%%%%%%%%%%%%%%%%%%%%%%%%%%%%%%%%%%%%%%%%%%%%%%%%%%%%%%%%%%%%%%%%%%%%%%
\begin{section}{Numerical Computations} \label{simulations}

In this Section, we test the optimization formulations of Sections \ref{energy_function_optimization}, \ref{MICP} on a model of natural gas network. Fig.~\ref{fig:toy_network} illustrates this exemplary network. There are $16$ nodes and $18$ edges. The nodes marked with green are sources and those with red are sinks. The arrows represent edges and the direction of the arrow denotes the nominal orientation of the edge. The cost $c_i$ is set to zero for sources and non-consumers. For consumers, we set $c_i$ proportional to the nominal demand at node $i$. We consider three different sets of upper bounds for the pressure squared variables $\pi$: $a)$ $\overline{\pi}_i = 10$ $b)$ $\overline{\pi}_i = 5$ $c)$ $\overline{\pi}_i = 3.5$. In all cases we set $\underline{\pi}_i = 0.5$. The slack node $0$ is set to the maximum pressure.
	
We solve the maximum throughput problem using the energy function formulation (\ref{reform_objective})-(\ref{reform_box}) and the MIQP formulation (\ref{MIQC_cost})-(\ref{MIQC_end}) both with and without compressors included into the set of optimization variables. To test the quality of the solution and bounds thus obtained, we compare the result with the one obtained with the help of publicly available global optimization software.
	
We implement the energy function formulation in IPOPT \cite{ipopt}, which is an interior point method for local optimization. We provide IPOPT with first derivative information by using Lemma \ref{duality_results}(a), and allow for automatic numerical Hessian computations. For the MIQP implementation, we use BONMIN \cite{bonmin}, which is a mixed integer convex programming software that uses IPOPT as the local search subroutine. We use SCIP \cite{scip}, a general purpose constraint optimization software, for global optimization.
	
Table \ref{tab:without_compression} shows the cost comparison between the energy function formulation, the MIQP relaxation and SCIP. The solution provided by the energy function formulation is a feasible solution and hence it is an upper bound. We conclude,  based on the experiments, that both the lower bound obtained by the MIQP and the upper bound from the energy function formulation are fairly close to the global optimum.
%	\textbf{Might add alternating minimization biconvex heuristics here}

	\begin{table}
	\caption{Numerical performance for the test without compression (See text for details.)}
	\label{tab:without_compression}
	\begin{center}
	\begin{tabular}{|c|c|c|c|}
		\hline
		\textbf{Pressure Bounds} $\rightarrow$& $0.5$-$5.0$ & $0.5$-$4.25$ & $0.5$-$3.5$ \\		
		\hline
		Energy function heuristics & -651 & -576 & -491 \\
		\hline
		Global Solver & -665 & -579 & -501 \\
		\hline
		MIQP & -685 & -581 & -501 \\
		\hline
	\end{tabular}
	\end{center}
	\end{table}
	
Table \ref{tab:with_compression} shows the optimal cost achieved by the energy function formulation, the MIQP relaxation and by SCIP.  Similar to the case without compression, one observes that the lower and upper bounds found by the MIQP and by the energy function based formulation are close to the global optimum.
		
	\begin{table}
	\caption{Numerical performance for the test with compression (See text for details.)}
	\label{tab:with_compression}
	\begin{center}
	\begin{tabular}{|c|c|c|c|}
		\hline
		\textbf{Pressure Bounds} $\rightarrow$& $0.5$-$5.0$ & $0.5$-$4.25$ & $0.5$-$3.5$ \\		
		\hline
		Energy function heuristics & -694 & -627 & -584 \\
		\hline
		Global Solver & -737  & -632 & -598 \\
		\hline
		MIQP & -752 & -640 & -598 \\
		\hline
	\end{tabular}
	\end{center}
	\end{table}

	\begin{figure}
	\centering
		\includegraphics[width=0.49\textwidth]{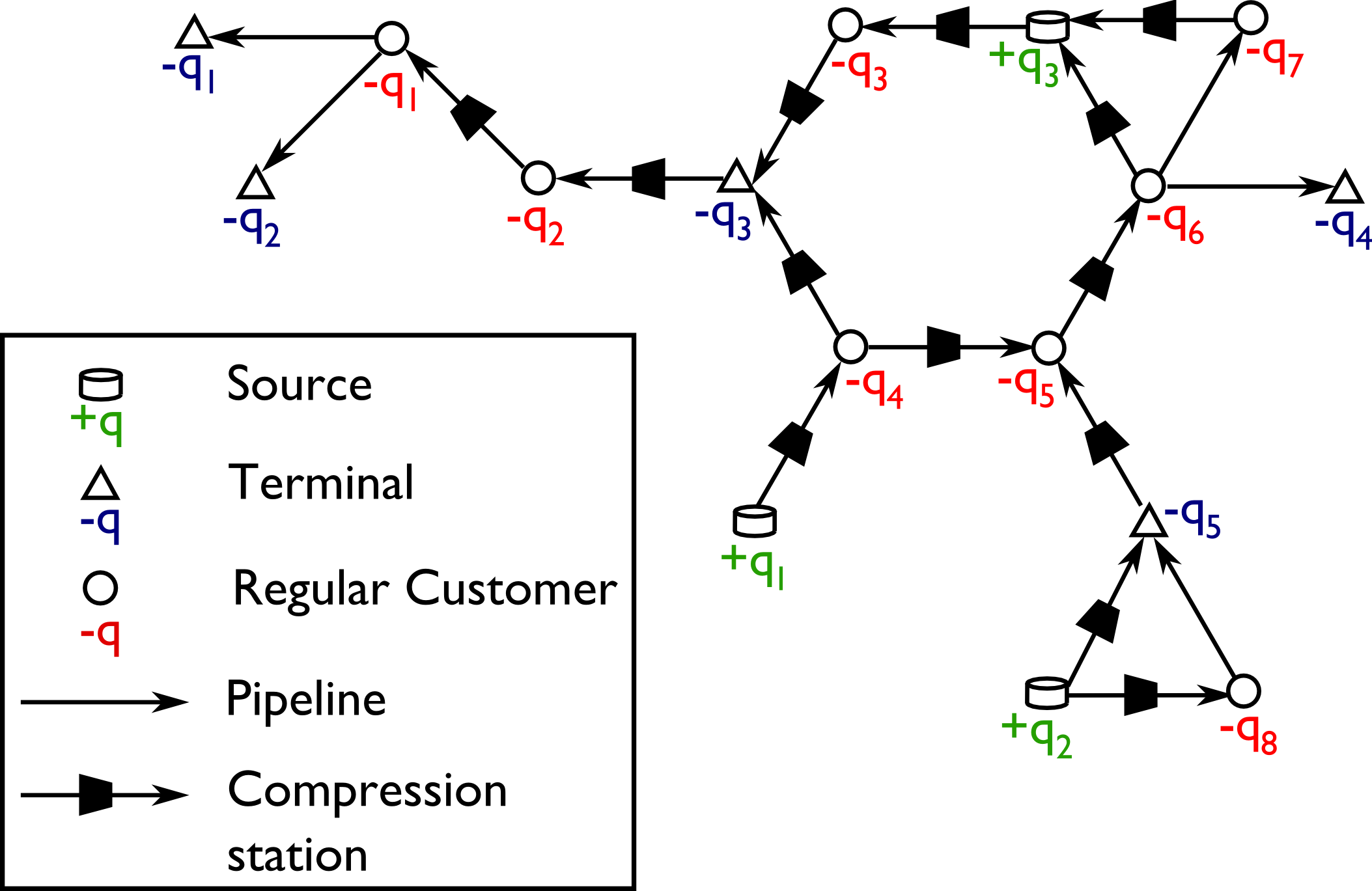}
		\caption{An exemplary (synthetic) gas network model}
	\label{fig:toy_network}
	\end{figure}	
\end{section}

%%%%%%%%%%%%%%%%%%%%%%%%%%%%%%%%%%%%%%%%%%%%%%%%%%%%%%%%%%%%%%%%%%%%%%%%%%%%%%%%
\section{Conclusions and Path Forward}
\label{sec:conclusions}
We developed a strictly convex variational representation of dissipative network flow equations that proves the uniqueness of the solution and provides an efficient algorithm to find it numerically. For the maximum throughput problem, we presented a twice differentiable optimization formulation using energy functions as a method to obtain a solution, and a mixed integer convex relaxation to test the optimality (or degree of sub-optimality) of the solution. We test both methods on a model of natural gas networks and observe that they were successful in producing solutions close to global optimality.

For future work, it would be useful to extend these techniques to dynamic optimization problems on dissipative flow networks (such as transient gas networks). The optimization formulations could also be used as subroutines in multi-level programs associated with network planning and expansion. Coupled with the monotonicity property, it also provides scope for extension to network optimization problems accounting for uncertainty such as the robust maximum throughput problem \cite{Vuffray15} and the problem of finding the most probable fluctuations that cause violation of bounds known as the instanton problem (see e.g. a similar setting discussed in the context of power systems \cite{11CPS}).

%{\color{red} ... for conclusions ... re-state theoretical results, remind the logic behind the relaxation,  mention success of the validation against numerical experiments ... for path forward ... extension of this technique to (a) dynamic problems, e.g. Model Predictive Control and Control formmulations; (b) bi- and more generally multi-level problems, including network design and planning; ... exploring other applications, such as traffic flows; (c) other unusual optimization problems, such as related to finding instantons (most probable failure configurations) ... }

%%%%%%%%%%%%%%%%%%%%%%%%%%%%%%%%%%%%%%%%%%%%%%%%%%%%%%%%%%%%%%%%%%%%%%%%%
\section{Acknowledgements}
%%%%%%%%%%%%%%%%%%%%%%%%%%%%%%%%%%%%%%%%%%%%%%%%%%%%%%%%%%%%%%%%%%%%%%%%%

The authors thank S. Backhaus for multiple discussions and advice, and A. Zlotnik for providing examplary model of the gas network he designed. The work at LANL was carried out under the auspices of the National Nuclear Security Administration of the U.S. Department of Energy at Los Alamos National Laboratory under Contract No. DE-AC52-06NA25396 and it was partially supported by DTRA Basic Research Project $\#10027-13399$. The authors also acknowledge partial support of the Advanced Grid Modeling Program in the US Department of Energy Office of Electricity.

\appendix{\emph{Proof of Monotonicity Property (Corollary \ref{monotonicity_property}):}
By Lemma \ref{duality_results}(a), $\nabla_{x} E^*(x, b) = \pi^*$, which means $\nabla_x \pi^* = \nabla_{xx}^2 E^*(x, b)$.
By Lemma \ref{duality_results}(b), we have the $\nabla_{xx}^2 E^*{x, b} = \left( \nabla_{\pi \pi}^2 E(\pi, b) \right)^{-1}$.
From (\ref{hessian_off_diag}) and (\ref{hessian_diag}) we observe that the Hessian $\nabla_{\pi \pi}^2 E(\pi, b)$ is a diagonally dominant invertible matrix with positive diagonal entries and non-positive off-diagonal entries. From a standard result in linear algebra we have that all entrees  in $\nabla_{xx}^2 E^*(x, b)$ are non-negative.

\bibliographystyle{IEEEtran}
\bibliography{Bib/GasFlow,Bib/Russian,Bib/RefConrado,Bib/GasFlowSM,Bib/cdc_gas_network}

\begin{thebibliography}{10}
\providecommand{\url}[1]{#1}
\csname url@samestyle\endcsname
\providecommand{\newblock}{\relax}
\providecommand{\bibinfo}[2]{#2}
\providecommand{\BIBentrySTDinterwordspacing}{\spaceskip=0pt\relax}
\providecommand{\BIBentryALTinterwordstretchfactor}{4}
\providecommand{\BIBentryALTinterwordspacing}{\spaceskip=\fontdimen2\font plus
\BIBentryALTinterwordstretchfactor\fontdimen3\font minus
  \fontdimen4\font\relax}
\providecommand{\BIBforeignlanguage}[2]{{%
\expandafter\ifx\csname l@#1\endcsname\relax
\typeout{** WARNING: IEEEtran.bst: No hyphenation pattern has been}%
\typeout{** loaded for the language `#1'. Using the pattern for}%
\typeout{** the default language instead.}%
\else
\language=\csname l@#1\endcsname
\fi
#2}}
\providecommand{\BIBdecl}{\relax}
\BIBdecl

\bibitem{77Nau}
J.~J. Maugis, ``Etude de réseaux de transport et de distribution de fluide [in
  french],'' pp. 243–--248, 1977.

\bibitem{12BNV}
\BIBentryALTinterwordspacing
F.~Babonneau, Y.~Nesterov, and J.-P. Vial, ``Design and operations of gas
  transmission networks,'' \emph{Operations Research}, 2012. [Online].
  Available:
  \url{http://or.journal.informs.org/content/early/2012/02/10/opre.1110.1001.abstract}
\BIBentrySTDinterwordspacing

\bibitem{Vuffray15}
M.~Vuffray, S.~Misra, and M.~Chertkov, ``Monotonicity of dissipative flow
  networks renders robust maximum profit problem tractable: General analysis
  and application to natural gas flows,'' in \emph{Submitted to 54th IEEE
  Conference on Decision and Control, Osaka}, 2015.

\bibitem{98Bol}
B.~Bollobas, \emph{Modern Graph Theory}.\hskip 1em plus 0.5em minus 0.4em\relax
  Springer, 1998.

\bibitem{58Ayl}
P.~Aylett, ``The energy-integral criterion of transient stability limits of
  power systems,'' \emph{Proceedings of the IEE-Part C: Monographs}, vol. 105,
  no.~8, pp. 527–--536, 1958.

\bibitem{81ASV}
A.~Araposthatis, S.~Sastry, and P.~Varaiya, ``Analysis of power-flow
  equation,'' \emph{International Journal of Electrical Power \& Energy
  Systems}, vol.~3, no.~3, pp. 115–--126, 1981.

\bibitem{81BH}
A.~Bergen and D.~Hill, ``A structure preserving model for power system
  stability analysis,'' \emph{Power Apparatus and Systems, IEEE Transactions
  on}, vol. 100, no.~1, pp. 25–--35, 1981.

\bibitem{15DLS}
K.~{Dvijotham}, S.~{Low}, and M.~{Chertkov}, ``{Convexity of Energy-Like
  Functions: Theoretical Results and Applications to Power System
  Operations},'' \emph{arxiv:1501.04052}, 2015.

\bibitem{Osiadacz_1987}
A.~Osiadacz, \emph{Simulation and analysis of gas networks}.\hskip 1em plus
  0.5em minus 0.4em\relax Gulf Publishing Company,Houston, TX, Jan 1987.

\bibitem{Wu2000197}
\BIBentryALTinterwordspacing
S.~Wu, R.~Ríos-Mercado, E.~Boyd, and L.~Scott, ``Model relaxations for the fuel
  cost minimization of steady-state gas pipeline networks,'' \emph{Mathematical
  and Computer Modelling}, vol.~31, no. 2–3, pp. 197 -- 220, 2000. [Online].
  Available:
  \url{http://www.sciencedirect.com/science/article/pii/S0895717799002320}
\BIBentrySTDinterwordspacing

\bibitem{Misra2015}
S.~Misra, M.~Fisher, S.~Backhaus, R.~Bent, M.~Chertkov, and F.~Pan, ``Optimal
  compression in natural gas networks: A geometric programming approach,''
  \emph{Control of Network Systems, IEEE Transactions on}, vol.~2, no.~1, pp.
  47--56, March 2015.

\bibitem{10CSADF}
G.~Como, K.~Savla, D.~Acemoglu, M.~A. Dahleh, and E.~Frazzoli, ``On robustness
  analysis of large-scale transportation networks,'' in \emph{Proceedings of
  the International Symposium on Mathematical Theory of Networks and Systems},
  2010, pp. 2399–--2406.

\bibitem{13Var_a}
P.~Varaiya, ``Max pressure control of a network of signalized intersections,''
  \emph{Transportation Research Part C: Emerging Technologies}, vol.~36, pp.
  177--195, 2013.

\bibitem{13Var_b}
------, ``The max-pressure controller for arbitrary networks of signalized
  intersections,'' in \emph{Advances in Dynamic Network Modeling in Complex
  Transportation Systems}.\hskip 1em plus 0.5em minus 0.4em\relax Springer,
  2013, pp. 27--66.

\bibitem{BoydCO}
S.~Boyd and L.~Vandenberghe, \emph{Convex Optimization}.\hskip 1em plus 0.5em
  minus 0.4em\relax Cambridge University Press, 2004.

\bibitem{Rockafellar70}
R.~T. Rockafellar, \emph{Convex Analysis}.\hskip 1em plus 0.5em minus
  0.4em\relax Princeton University Press, 1970.

\bibitem{Merris94}
R.~Merris, ``Laplacian matrices of graphs: a survey,'' \emph{Linear Algebra
  Applications}, vol. 197, 1994.

\bibitem{Chaiken78}
S.~Chaiken and D.~J. Kleitman, ``Matrix tree theorems,'' \emph{Journa of
  Combinatorial Theory}, vol.~24, no. 377-381, 1978.

\bibitem{McCormick76}
G.~McCormick, ``Computability of global solutions to factorable nonconvex
  programs: {P}art {I} - convex underestimating problems,'' \emph{Mathematical
  Programming}, vol.~10, no. 146-175, 1976.

\bibitem{68WL}
P.~Wong and R.~Larson, ``Optimization of natural-gas pipeline systems via
  dynamic programming,'' \emph{Automatic Control, IEEE Transactions on},
  vol.~13, no.~5, pp. 475--481, 1968.

\bibitem{10Bor}
C.~Borraz-S\'{a}nchez, ``Optimization methods for pipeline transportation of
  natural gas,'' Ph.D. dissertation, Department of Informatics, University of
  Bergen, Norway, October 2010.

\bibitem{Misra13}
S.~Misra, M.~W. Fisher, R.~Bent, S.~Backhaus, M.~Chertkov, and F.~Pan,
  ``Optimal compression in natural gas networks: a geometric programming
  approach,'' \emph{IEEE Control of Network Systems}, 2014.

\bibitem{ipopt}
A.~W. L.~T. Biegler, ``On the implementation of a primal-dual interior point
  filter line search algorithm for large-scale nonlinear programming,''
  \emph{Mathematical Programming}, vol. 106, no.~1, pp. 25--27, 2006.

\bibitem{bonmin}
P.~Bonami, A.~Wachter, L.~T. Biegler, A.~R. Conn, G.~Cornuejols, I.~E.
  Grossman, C.~D. Laird, J.~Lee, A.~Lodi, F.~Margot, and N.~Sawaya, ``An
  algorithmic framework for convex mixed integer nonlinear programs.''
  \emph{IBM Research Report, RC23771}, 2005.

\bibitem{scip}
T.~Achterberg, T.~Berthold, T.~Koch, and K.~Wolter, ``Constraint integer
  programming: a new approach to integrate {CP} and {MIP},'' \emph{Integration
  of AI and OR techniques in constraint programming for combinatorial
  optimization problems, CPAIOP}, 2008.

\bibitem{11CPS}
M.~Chertkov, F.~Pan, and M.~Stepanov, ``Predicting failures in power grids: The
  case of static overloads,'' \emph{Smart Grid, IEEE Transactions on}, vol.~2,
  no.~1, pp. 162--172, March 2011.

\end{thebibliography}

%%%%%%%%%%%%%%%%%%%%%%%%%%%%%%%%%%%%%%%%%%%%%%%%%%%%%%%%%%%%%%%%%%%%%%%%%%%%%%%%

%\bibliographystyle{IEEEtran}
%\bibliography{20150223_bib_LR}

%References are important to the reader; therefore, each citation must be complete and correct. If at all possible, references should be commonly available publications.
%
%\begin{thebibliography}{99}
%
%\bibitem{c1}
%J.G.F. Francis, The QR Transformation I, {\it Comput. J.}, vol. 4, 1961, pp 265-271.
%
%\bibitem{c2}
%H. Kwakernaak and R. Sivan, {\it Modern Signals and Systems}, Prentice Hall, Englewood Cliffs, NJ; 1991.
%
%\bibitem{c3}
%D. Boley and R. Maier, "A Parallel QR Algorithm for the Non-Symmetric Eigenvalue Algorithm", {\it in Third SIAM Conference on Applied Linear Algebra}, Madison, WI, 1988, pp. A20.
%
%\end{thebibliography}

\end{document}